\documentclass[11pt]{article}
\usepackage{cite}
\usepackage{mathrsfs}
\usepackage{amsfonts}
\usepackage{amsmath}
\usepackage{amsfonts,amssymb,color}
\usepackage{dsfont}
\usepackage{curves}
\usepackage{mathrsfs}
\usepackage{pifont}
\usepackage{amssymb}
\usepackage{latexsym,amsmath,amssymb,amsfonts,epsfig,graphicx,cite,psfrag}
\usepackage{eepic,color,colordvi,amscd}

\newtheorem{theorem}{Theorem}[section]

\newtheorem{lemma}{Lemma}[section]

\newtheorem{claim}{Claim}[section]
\newtheorem{conjecture}{Conjecture}[section]

\newcommand{\qed}{\hfill\rule{0.5em}{0.809em}}

\def\qed{\hfill \rule{4pt}{7pt}}

\def\pf{\noindent {\it Proof. }}

\begin{document}

 \title{Borodin-Kostochka conjecture for a family of $P_6$-free graphs}
  \author{Di Wu$^{1,}$\footnote{Email: 1975335772@qq.com},  \;\; Rong Wu$^{2,}$\footnote{Corresponding author: wurong@sjtu.edu.cn}\\\\
 	\small $^1$Department of Mathematics and Physics\\
 	\small Nanjing Institute of Technology, 1 Hongjing Avenue, Nanjing, 211167, China\\
    \small $^2$School of Mathematical Sciences\\
    \small Shanghai Jiao Tong University, 800 Dongchuan Road, Shanghai 200240, China}
 \date{}

 \maketitle
\begin{abstract}
Borodin and Kostochka conjectured that every graph $G$ with $\Delta\ge9$ satisfies $\chi\le$ max $\{\omega, \Delta-1\}$. Gupta and Pradhan proved the Borodin-Kostochka
conjecture for ($P_5$, $C_4$)-free graphs [{\em J. Appl. Math. Comp.} \textbf{65} (2021) 877-884]. 
In this paper, we prove the Borodin-Kostochka
conjecture for ($P_6$, apple, torch)-free graphs, that is, graphs
with no induced $P_6$, no induced $C_5$ with a hanging edge, and no induced $C_5$ and $C_4$ sharing exactly an induced $P_3$. This generalizes the result of Gupta and Pradhan from the perspective of allowing the existence of $P_5$.
\begin{flushleft}
\end{flushleft}

\end{abstract}

\section{Introduction}

This paper focuses on finite and simple graphs. We denote a {\em path} and a {\em cycle} on $k$ vertices by $P_k$ and $C_k$, respectively, and refer to \cite{BM08} for any notations and terminology that are not defined here. 

Let  $G$ be a graph, and let  $X$ be a subset of $V(G)$. We denote by  $G[X]$ the subgraph of  $G$ that is induced by $X$. We refer to  $X$ as a {\em clique} if  $G[X]$ forms a complete graph, and as a {\em stable set} if  $G[X]$ contains no edges. The {\em clique number} $\omega(G)$ of $G$ is defined as the largest size of any clique within  $G$.

For any vertex $v\in V(G)$, let $N_G(v)$ represent the set of vertices that are adjacent to $v$, and let $d_G(v)$ represent the number of vertices in $N_G(v)$. If there is no potential confusion, we will exclude the subscript $G$ and use the notation $N(v)$ instead. We denote the maximum degree of $G$ by  $\Delta(G)$.
 
Consider two graphs, $G$ and $H$, that are vertex-disjoint from each other. The {\em union} of  $G $ and $H$, denoted $G\cup H$, is defined as the graph whose vertex set is $V(G)\cup V(H)$ and whose edge set is $E(G)\cup E(H)$. We say that $G$ induces $H$ when there exists an induced subgraph of $G$ that is isomorphic to $H$. Conversely, we state that $G$ is $H$-free if it does not induce $H$. Similarly, for a family of graphs $\mathcal{H}$, we define $G$ to be $\mathcal{H}$-free if it does not induce any graph from the family $\mathcal{H}$.

For vertices $u, v\in V(G)$, we denote $u\sim v$ to indicate that $uv\in E(G)$ and write $u\not\sim v$ when $uv\not\in E(G)$. A {\em hole} in a graph $G$ refers to an induced cycle with a length of at least 4. A {\em $k$-hole} is defined as a hole that specifically has a length of $k$. If $k$ is odd, the $k$-hole is termed an {\em odd hole}; if $k$ is even, it is referred to as an {\em even hole}.

A {\em $k$-coloring} of $G$ is a partition of $V(G)$ into $k$ stable sets, and the {\em chromatic number} $\chi(G)$ of $G$ is the minimum integer $k$ such that $G$ admits a $k$-coloring.
A simple lower bound for the chromatic number $\chi(G)$ is $\omega(G)$. Additionally, by using a greedy coloring method, a straightforward upper bound for $\chi(G)$ is $\Delta(G)+1$. In 1941, Brooks \cite{B41} noted that for a graph $G$, the chromatic number is at most $\Delta(G)$, except when $G$ is a complete graph or an odd cycle. It is clear that if $G$ is a complete graph or an odd cycle, then $\chi(G)=\Delta(G)+1$.

\begin{theorem}\label{Brook}\cite{B41}
	Let $G$ be a graph with $\Delta(G)\geq 3$. Then $\chi(G)\leq $ max$\{\Delta(G),\omega(G)\}$.
\end{theorem}

In 1977, Borodin and Kostochka \cite{BK77} conjectured that Brooks' bound could be further enhanced when $ \Delta(G) \geq 9 $.

\begin{conjecture}\label{BKC}\cite{BK77}
	Let $G$ be a graph with $\Delta(G)\geq 9$. Then $\chi(G)\leq $ max$\{\Delta(G)-1,\omega(G)\}$.
\end{conjecture}

According to Brooks' Theorem, any graph $G$ with $\chi(G) > \Delta(G) \geq 9$ must contain $ K_{\Delta(G)+1}$. Therefore, the Borodin-Kostochka conjecture is equivalent to the assertion that every graph $G$ with $\chi(G) = \Delta(G) \geq 9$ includes $K_{\Delta(G)}$. This conjecture seems to be challenging, as its result is so strong that it was demonstrated in \cite{CLR22} that the Borodin-Kostochka conjecture cannot be enhanced by imposing $\Delta(G)\geq 8$ or $\omega(G)\le \Delta(G)-2$. The current best result on this conjecture is that Reed \cite{R99} demonstrated this conjecture for all graphs $G$ where $\Delta(G)$ is sufficiently large, in 1999.

\begin{theorem} \label{RS}\cite{R99}
Every graph with $\chi(G)=\Delta(G)\geq 10^{14}$ contains $K_{\Delta(G)}$.
\end{theorem}

\begin{figure}[htbp]\label{fig-1}
	\begin{center}
		\includegraphics[width=12cm]{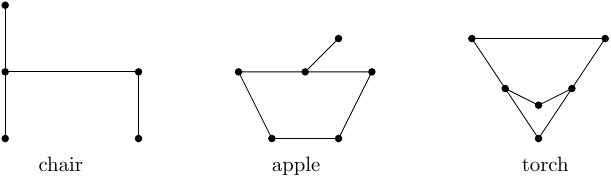}
	\end{center}
	\vskip -15pt
	\caption{Illustration of chair, apple and torch.}
\end{figure}

In 2013, Cranston and Rabern \cite{CR13} proved the Borodin-Kostochka conjecture for $K_{1,3} $-free graphs. Recently, together with Lafayette \cite{CLR22}, they extended the result to ($P_5$, gem)-free graphs. Additionally, Gupta and Pradhan \cite{GP21} demonstrated the conjecture for $(P_5, C_4)$-free graphs, while Chen, Lan, Lin and Zhou \cite{cllz} proved it for odd-hole-free graphs. Wang \cite{W23} established the conjecture for some $(P_2 \cup P_3)$-free graphs. Furthermore, Dhurandhar \cite{D} proved it for forbidden subgraphs with five vertices, specifically for $P_4 \cup K_1 $, $P_5$, and the chair graph (which is derived from $K_{1,3}$ by subdividing an edge once, see Figure \ref{fig-1} for chair).

In this paper, we want to focus on forbidden subgraphs with six vertices that prove the Borodin-Kostochka
conjecture for ($P_6$, apple, torch)-free graphs, where apple is an induced $C_5$ with a hanging edge and torch is an induced $C_5$ and an induced $C_4$ sharing exactly an induced $P_3$ (See Figure \ref{fig-1} for apple and torch). The result is stated in the following theorem.

\begin{theorem} \label{mian theorem}
Let $G$ be a $($$P_6$, apple, torch$)$-free graph. If $\Delta(G)\geq 9$, then $\chi(G)\leq $ max$\{\Delta(G)-1,\omega(G)\}$.
\end{theorem}

Obviously, the class of ($P_6$, apple, torch)-free graphs is a superclass of ($P_5, C_4$)-free graphs and allows the existence of $P_5$. Therefore, this generalizes the result of Gupta and Pradhan in \cite{GP21} from the perspective of allowing the existence of $P_5$.

\section{Notations and preliminary results}

A graph class $\cal{A}$ is called {\em hereditary} if for every graph $G \in \cal{A}$, every induced subgraph $H$ of $G$ also belongs to $\cal{A}$. 
Every class of graphs characterized by a list of forbidden induced subgraphs is a hereditary class; 
in particular, the class of ($P_6$, apple, torch)-free graphs is hereditary. 

For a graph $G$,  if $\chi(H)<\chi(G)$ for every proper induced subgraph $H$ of $G$, then we say $G$ is {\em vertex-critical}. 

In our proof, we will rely on the following two lemmas, which provide essential support that will be referenced and utilized throughout our argumentation.

\begin{lemma}{\em\cite{C76,CLR22,K80}}\label{kc}
Let $\cal{G}$ be a hereditary class of graphs. If the
Borodin-Kostochka conjecture is false for some $G\in \cal{G}$, then it is false for some $G\in \cal{G}$ with $\Delta(G)=9$.
\end{lemma}

\begin{lemma}{\em \cite{GP21}}\label{GP}
	If $G$ is a minimal counterexample for the Borodin-Kostochka conjecture with $\Delta(G)=9$, then $G$ is vertex-critical.
\end{lemma}

According to Lemma \ref{kc}, we can always assume that \(\Delta(G) = 9\) for every graph \(G\) discussed in this paper. Furthermore, by Lemma \ref{GP}, we can consider each minimal counterexample to the Borodin-Kostochka conjecture with \(\Delta(G) = 9\) as a vertex-critical graph.

\medskip

Let $C$ be some coloring of $G$.  If $u\in V(G)$ is colored $m$ in $C$, then $u$ is called an {\em $m$-vertex}.
If $v\in N(u)$ is a $r$-vertex, then we call $v$ is a $r$-vertex of $u$. Particularly, if $v$ is the unique $r$-vertex in $N(u)$,
then $r$ is called a {\em unique color} of $u$, $v$ is called a unique $r$-vertex of $u$ and if $N(u)$ has more than one $r$-vertex,
then $r$ is called a {\em repeat color} of $u$.
If the number of colors in the coloring of $N(u)$ is smaller than the number of colors in $C$, we refer to $u$ as having a {\em missing color}. 

Furthermore, for two distinct colors $i$ and $j$ in $C$, and two distinct vertices $x$ and $y$ in $V(G)$, an $xy$-path is referred to as an $(i,j)$-$xy$-{\em path} if all its vertices are colored alternately with $i$ and $j$.

Recall that $\Delta(G)=9$. We select a vertex $u\in V(G)$ with a degree of 9 and let $N(u)=\{u_1,u_2,\dots,u_7,x,y\}$. Let $\phi(V(G)\setminus\{u\})=\{1,2,\dots,8\}$ be an 8-coloring of $G-u$, where we assign color $i$ to vertex $u_i$ for $i\in\{1,2,\dots,7\}$, and we assign color 8 to vertices $x$ and $y$. If there exists a vertex $u$ and a coloring $\phi$ satisfying the conditions mentioned above, we say that $G$ has a $(u,\phi)$.

\begin{lemma}\label{relax}
If $G$ is a minimal counterexample for the Borodin-Kostochka conjecture with $\Delta(G)=9$, then $G$ must have a $(u,\phi)$ and satisfy the following three properties. 

\medskip

{\textbf{\em(a)}} For $1\le i\le7$, $u_i$ cannot have missing color.

\medskip

{\textbf{\em(b)}} For $1\le i,j\le7$ and $i\ne j$, let $u_i\not\sim u_j$. Then $G$ must have an $(i,j)$-$u_iu_j$-path. Specifically, if $P$ is such a path, then its length is odd and the induced subgraph $G[V(P)\cup\{u\}]$ forms an odd hole.

\medskip

{\textbf{\em(c)}} For $1\le i,j\le7$ and $i\ne j$, let $u_i\not\sim u_j$. Then $u_i$ has at least one $j$-vertex such that the $j$-vertex has no missing color and has an $i$-vertex other than $u_i$.
\end{lemma}
 
\medskip

\pf According to Lemma \ref{kc}, it is evident that there exists a minimal counterexample $G$ for the Borodin-Kostochka conjecture with $\Delta(G)=9$. 
Choose $u$ with $d(u)=9$. Since $G$ is a minimal counterexample for the Borodin-Kostochka conjecture, $\chi(G)=9$ and $\chi(G-u)=8$. Therefore, $G$ has a $(u,\phi)$.

Suppose {\textbf{(a)}} is false, indicating that $u_1$ has a missing color $r$. In this case, we assign the color $r$ to $u_1$ and the color $1$ to $u$. Consequently, $G$ becomes $8$-colorable, which contradicts the assumption. So, {\textbf{(a)}} holds.

Suppose {\textbf{(b)}} is not true. Without loss of generality, let $i=1,j=2$. Suppose $G$ has no $(1,2)$-$u_1u_2$-path. Let $S$ be the set of vertices in $G$ that are colored 1 or 2. Let $T$ be the component of $G[S]$ that contains the vertex $u_1$. Since $G$ has no $(1,2)$-$u_1u_2$-path, we can conclude that $u_2$ is not contained in $T$. If we interchange the colors of the vertices in $T$, specifically coloring $u_1$ with color 2, we obtain a new coloring of $G$. However, this new coloring implies that $G$ is $8$-colorable, which contradicts our assumption that $G$ is not $8$-colorable. Therefore,  {\textbf{(b)}} is true.

If {\textbf{(c)}} is false, we can assume without loss of generality that $i = 1$ and $j = 2$. By {\textbf{(a)}}, $u_1$ has at least one $2$-vertex. If $u_1$ has more than two $2$-vertices, then the total number of neighbors of $u_1$ would exceed 8, which is not possible since $d_{G-u}(u_1) \leq 8$. Hence, $u_1$ can have at most two $2$-vertices.

First, let us assume that $u_1$ has a unique 2-vertex, which we will denote as $v(2)$. If $v(2)$ has a missing color $r$, we can color $v(2)$ with $r$, color $u_1$ with 2, and color $u$ with 1, which results in a contradiction. Therefore, we conclude that $v(2)$ has no missing color. If $v(2)$ has the unique 1-vertex $u_1$, then we can exchange the colors of $v(2)$ and $u_1$, while coloring $u$ with 1, which leads to a contradiction. Therefore, we conclude that $v(2)$ must have a 1-vertex other than $u_1$.

Then, we suppose that $u_1$ has two 2-vertices, denoted as $v(2)$ and $v'(2)$. It is clear that $\{v(2),v'(2)\}$ forms a stable set. If both $v(2)$ and $v'(2)$ have the unique 1-vertex $u_1$, we can color  $v(2)$ and $v'(2)$ with 1, color $u_1$ with 2, and color $u$ with 1, which leads to a contradiction. Therefore, at least one of the elements in $\{v(2),v'(2)\}$ must have a 1-vertex that is distinct from $u_1$.

If $v(2)$ has a missing color $r$ and $v'(2)$ has a missing color $r'$ , then  we can assign color $r$ to $v(2)$, color $r'$ to $v'(2)$, color 2 to $u_1$, and color 1 to $u$. This leads to a contradiction. Therefore, at least one element of $\{v(2),v'(2)\}$ has no missing color.

Now, if $v(2)$ or $v'(2)$ has no missing color and has an $i$-vertex other than $u_i$, then {\textbf{(c)}} is true. 
So, we may assume that $v(2)$ has a missing color $r$ and has a $1$-vertex other than $u_1$, but $v'(2)$ has no missing color and has the unique $1$-vertex $u_1$. Then, we color $v(2)$ by $r$, color $v'(2)$ by 1, color $u_1$ by 2, and color $u$ by 1. This leads to a contradiction. Therefore,  {\textbf{(c)}} holds.
\qed

\medskip

By Lemma \ref{relax}, we can conclude that every minimal counterexample $G$ for the Borodin-Kostochka conjecture with $\Delta(G)=9$ has a $(u,\phi)$.  We refer to such $G$ a {\em relaxed graph} with a $(u,\phi)$, in particular, every relaxed graph satisfies properties {\textbf{(a)}},{\textbf{(b)}} and {\textbf{(c)}}.

\medskip

The following theorem was proven by Dhurandhar in \cite{D}, but the proof contains some minor flaws. To make it easier for readers to understand, we provide a more comprehensible proof.

\begin{theorem}\label{MD}\cite{D}
	Let $G$ be a relaxed graph with a $(u,\phi)$. Then at least one of the following statements does not hold.
	
	\medskip
	
	\textbf{(\romannumeral 1)} $u_i$ is adjacent to at least four $u_k$'s for $1\le i,k\le7$.
	
	\medskip
	
	\textbf{(\romannumeral 2)} For $1\le i\le7$, $u_i\sim x$ or $u_i\sim y$.
\end{theorem}
\pf Suppose, to the contrary, that $G$ satisfies both statements \textbf{(\romannumeral 1)} and \textbf{(\romannumeral 2)}. Since $G$ is a minimal counterexample for the Borodin-Kostochka conjecture, we have that $\chi(G)=9, \omega(G)\le8$ and $\chi(G-u)=8$. To prove the statement, we divide the proof process into two cases based on whether $\{u_1,u_2,\dots,u_7\}$ forms a clique or not.

\medskip

\noindent{\bf Case} 1. $\{u_1,u_2,\dots,u_7\}$ is not a clique.

Without loss of generality, we may assume that \( u_1 \not\sim u_2 \). Therefore, we will proceed to prove that

\begin{equation}\label{eqa-1}
\mbox{for $i\in\{3,4,\dots,7\}$, $u_i$ cannot be the unique $i$-vertex of both $u_1$ and $u_2$}.
\end{equation}

Suppose otherwise; let us assume that $u_3$ is the unique $3$-vertex of both $u_1$ and $u_2$. According to \textbf{(a)}, $u_3$ cannot have two repeat colors, and thus we suppose $u_1$ is the unique 1-vertex of $u_3$. Now, we can color $u_1,u_2$ by 3, $u_3$ by 1, $u$ by 2, which leads to a contradiction. This establishes the proof for (\ref{eqa-1}).

Next, we will prove that 

\begin{equation}\label{eqa-2}
\mbox{$|N(u_1)\cap N(u_2)\cap N(u)|\le2$}.
\end{equation}

If this is not the case, we may assume that $\{v_1,v_2,v_3\}\subseteq N(u_1)\cap N(u_2)\cap N(u)$.
According to \textbf{(a)}, $u_1$ has at most one repeated color, and the same applies to $u_2$. Consequently, $u_1$ and $u_2$ must share at least one common unique vertex among $\{v_1,v_2,v_3\}$, which we can denote as $v_1$.  By (\ref{eqa-1}), $v_1\not\in\{v_3,v_4,\dots,v_7\}$, which implies that $v_1\in\{x,y\}$. By symmetry, we can assume that $v_1 = x$. If $x$ has a missing color $r$, then we could color $x$ with $r$, $u_1$ with 8, and $u$ with 1, leading to a contradiction. Thus, $x$ cannot have a missing color; in other words, $x$ has at most one repeated color. Without loss of generality, we can assume that $u_1$ is the unique 1-vertex of $x$. In this case, we could color $u_1$ and $u_2$ with color 8, $x$ with color 1, and $u$ with color 2. This also leads to a contradiction, thereby proving (\ref{eqa-2}).
 
By \textbf{(\romannumeral 1)} and given that \(d(u) = 9\), we conclude that $|N(u_1)\cap N(u_2)\cap N(u)|\ge3$. This contradicts the result established in (\ref{eqa-2}).

\medskip

\noindent{\bf Case} 2. $\{u_1,u_2,\dots,u_7\}$ is a clique.

By \textbf{(\romannumeral 2)}, each vertex in $\{u_1,u_2,\dots,u_7\}$ is adjacent to $x$ or $y$. So, we may assume that $|N(x)\cap\{u_1,u_2,\dots,u_7\}|\ge|N(y)\cap\{u_1,u_2,\dots,u_7\}|$, which implies that $|N(x)\cap\{u_1,u_2,\dots,u_7\}|\ge4$. 

Since $\omega(G)\le8$, there must be a vertex in $\{u_1,u_2,\dots,u_7\}$ that is not adjacent to $x$. Let us denote this vertex as $u_1$.  Additionally, we can conclude that $u_1\sim y$. As $|N(x)\cap\{u_1,u_2,\dots,u_7\}|\ge|N(y)\cap\{u_1,u_2,\dots,u_7\}|$, we have that there exists a vertex in $\{u_1,u_2,\dots,u_7\}$ which is not adjacent to $y$, say $u_2$.

If $x$ has a missing color, say color $r$, then we can color $x$ with color $r$. This leads us to Case 1, a contradiction. Hence, $x$ has no missing color. Given that $|N(x)\cap\{u_1,u_2,\dots,u_7\}|\ge4$, we can deduce that $x$ and $u_1$ share at least two unique vertices from $\{u_1,u_2,\dots,u_7\}$, which we can denote as $u_i$ and $u_j$. If $x$ is the unique 8-vertex of $u_i$, then we can assign color $i$ to $x$ and $u_1$, color $u_i$ with 8, and color $u$ with $1$. This leads a contradiction. Therefore, $u_i$ must have two 8-vertices. In this case, we can color $x$ and $u_1$ with color $i$, and $u_i$ with color 1. However, this would lead us back to Case 1 as $y\not\sim u_2$, which is a contradiction.

\medskip

This completes the proof of Theorem \ref{MD}.
\qed

\medskip

\section{Proof of Theorem \ref{mian theorem}}

In this section, we will provide a proof for Theorem \ref{mian theorem}.

\medskip

Let us begin by assuming the contrary and considering a minimal counterexample $G$ to the Borodin-Kostochka conjecture. Recall that $G$ is a relaxed graph with a $(u,\phi)$ which satisfies properties {\textbf{(a)}}, {\textbf{(b)}} and {\textbf{(c)}}. In order to prove Theorem~\ref{mian theorem},  by Theorem~\ref{MD}, we only need to establish Claim~\ref{claim-1} and Claim~\ref{claim-2}.

\begin{claim}\label{claim-1}
	$u_i$ is not adjacent to at most two $u_k$'s for $1\le i,k\le7$.
\end{claim}

\pf Suppose to its contrary. Based on symmetry, we may assume that $u_1\not\sim u_2$, $u_1\not\sim u_3$ and $u_1\not\sim u_4$. By \textbf{(c)}, $u_2$ has a 1-vertex $v(1)$ which has no missing color and has a $2$-vertex $w(2)$ other than $u_2$. As $G$ is $P_6$-free and apple-free, we may assume $w(2)\sim u_1$ and $\{u,u_2,v(1),w(2),u_1\}$ induces a $C_5$. Similarly, we observe that there is a $C_5$ induceded by $\{u,u_3,t(1),s(3),u_1\}$, as well as another $C_5$ induced by $\{u,u_4,q(1),p(4),u_1\}$.

If $v(1)=t(1)=q(1)$, then $v(1)$ would have two 2-vertices, two 3-vertices, and two 4-vertices. However, since $d_{G-u}(v(1))\leq 9$, we know that $v(1)$ must have a missing color, which contradicts the assumption. Therefore, we can assume that $v(1)\neq t(1)$.

To forbid an induced $P_6$ on $\{v(1),w(2),u_1,u,u_3,t(1)\}$, we can conclude that $v(1)\sim u_3$ or $w(2)\sim t(1)$ or $w(2)\sim u_3$. If $w(2)\sim t(1)$, then $u_2\sim t(1)$; otherwise, $\{v(1),u_2,u,u_1,w(2),t(1)\}$ would induce an apple, resulting in a contradiction. However, this means that the vertices $\{v(1),u_2,u,u_1,w(2),t(1)\}$ would now induce a torch, which is again contradictory. So, $w(2)\not\sim t(1)$. 

If $v(1)\sim u_3$, then $v(1)\sim s(3)$ as otherwise $\{v(1),u_3,u,u_1,s(3),t(1)\}$ induces an apple, which is a contradiction. However, this would also result in $\{v(1),u_3,u,u_1,s(3),t(1)\}$ inducing a torch, which is contradictory.

Therefore, we may always assume that $w(2)\sim u_3$. Consequently, to forbid a torch on $\{v(1),u_2,u,u_1,w(2),u_3\}$, we have that $u_2\sim u_3$. Similarly, we have that $u_2\sim s(3)$. Since $q(1)\ne v(1)$ or $q(1)\ne t(1)$, we may assume that $q(1)\ne v(1)$. Based on the previous observations, we can deduce that $u_2\sim u_4$ and $u_2\sim q(4)$. However, this configuration results in $u_2$ having two 3-vertices ($u_3$ and $s(3)$) and two 4-vertices ($u_4$ and $q(4)$), which implies that $u_2$ has a missing color, a contradiction. \qed

\begin{claim}\label{claim-2}
	For $1\le i\le7$, $u_i\sim x$ or $u_i\sim y$.
\end{claim}

\pf Suppose to its contrary. By symmetry, $u_1\not\sim x$ and $u_1\not\sim y$. By \textbf{(a)}, $u_1$ has a 8-vertex $v(8)$ which has a 1-vertex $w(1)$ other than $u_1$ as otherwise we can color $u_1$ by 8 and color $u$ by 1, a contradiction. If both $x$ and $y$ have no 1-vertex, we can color $x$ and $y$ by 1, and color $u$ by 8, a contradiction. So, by symmetry, we may assume that $x$ has a 1-vertex $v(1)$. Then $w(1)=v(1)$ or $w(1)\sim x$ as otherwise $\{w(1),v(8),u_1,u,x,v(1)\}$ induces a $P_6$, a contradiction. Therefore, no matter which case happens, $\{u,u_1,v(8),w(1),x\}$ induces a $C_5$. Now, $\{u,u_1,v(8),w(1),x,y\}$ induces an apple if $w(1)\not\sim y$ and $\{u,u_1,v(8),w(1),x,y\}$ induces a torch if $w(1)\sim y$, both are contradictions. \qed

\medskip

This completes the proof of Theorem \ref{mian theorem}.
\qed

\bigskip

\noindent{\bf Remarks}. Based on Lemma \ref{relax}\textbf{(b)} and Theorem \ref{MD}, it is straightforward to conclude that the Borodin-Kostochka conjecture holds for odd-hole-free graphs; because if it is not true, then the graph $G[V(P)\cup\{u\}]$ in Lemma \ref{relax}\textbf{(b)} must form an odd hole, leading to a contradiction.

A vertex $v$ in a graph $G$ is considered {\em bisimplicial} if its neighborhood $N(v)$ can be divided into two cliques. A graph is defined as {\em quasi-line} if every vertex within the graph is bisimplicial. According to reference \cite{CS12}, every quasi-line graph is $K_{1,3}$-free. Additionally, Cranston and Rabern \cite{CR13} proved the Borodin-Kostochka conjecture for $K_{1,3}$-free graphs,  and consequently for every quasi-line graph. Based on the reference \cite{CS23}, it is established that every even-hole-free graph contains a bisimplicial vertex. Therefore, 
it is reasonable to believe that the Borodin-Kostochka conjecture holds true for even-hole-free graphs. This is still open.

Lan, Liu, and Zhou \cite{LLZ24} established the Borodin-Kostochka conjecture for ($P_2 \cup P_3$, $C_4$)-free graphs, while Gupta and Pradhan \cite{GP21} proved the conjecture for ($P_5$, $C_4$)-free graphs. A natural question arises: does the Borodin-Kostochka conjecture hold for ($P_6$, $C_4$)-free graphs? This looks interesting.
\bigskip

\noindent{\bf Acknowledgement}: We thank Professor Zixia Song for helpful suggestions and for informing us about some problems related to even-hole-free graphs. We also  thank the anonymous referees for their helpful suggestions and comments that improved the presentation of the paper greatly.

\section*{Declarations}

\begin{itemize}
	\item \textbf{Funding}\quad Research of the first author was supported  by the Scientific Research
	Foundation of Nanjing Institute of Technology, China (No. YKJ202448), research of the second author was supported by National Key R\&D Program of China under Grant No.2022YFA1006400 and Shanghai Municipal Education Commission (No. 2024AIYB003).
	\item \textbf{Conflict of interest}\quad The authors declare no conflict of interest.
	\item \textbf{Data availibility statement}\quad This manuscript has no associated data.
\end{itemize}


\begin{thebibliography}{9999}
	
	\bibitem{BM08}  J. A. Bondy, U. S. R. Murty, Graph Theory, Springer, New York, 2008.
	
	
    \bibitem{BK77} O.V. Borodin, A. V. Kostochka, On an upper bound of a graph's chromatic number, depending on the
graph's degree and density, J. Comb. Theory. Ser. B, 23 (1977) 247-250.

    \bibitem{B41} R.L. Brooks, On colouring the nodes of a network, Math. Proc. Camb. Phil. Soc., 37 (1941) 194-197.

    \bibitem{C76} P.A. Catlin, Embedding subgraphs and coloring graphs under extremal degree conditions. Ann Arbor, MI,
ProQuest LLC, Ph.D. thesis, The Ohio State University, 1976.

    \bibitem{cllz} R. Chen, K. Lan, X. Lin, Y. Zhou, Borodin-Kostochka Conjecture holds for odd-hole-free graphs, Graphs Comb., 40 (2024) 26.
    
   \bibitem{CS12} M. Chudnovsky, P. Seymour, Claw-free graphs. VII. Quasi-line graphs, J. Comb. Theory. Ser. B, 102 (2012) 1267-1294.

   \bibitem{CS23} M. Chudnovsky, P. Seymour, Even-hole-free graphs still have bisimplicial vertices, J. Comb. Theory. Ser. B, 161 (2023) 331-381.
    
    \bibitem{CR13} D.W. Cranston, L. Rabern, Coloring claw-free graphs with $\Delta-1$ colors, Siam J. Disc. Math., 27 (2013) 534-549.

    \bibitem{CLR22} D.W. Cranston, H. Lafayette, L. Rabern, Coloring ($P_5$, gem)-free graphs with $\Delta-1$ colors, J. of Graph Theory, 100 (2022) 1-10.
    
    \bibitem{D} M. Dhurandhar, Validity of Borodin and Kostochka Conjecture for classes of graphs without a single, forbidden subgraph on 5 vertices. arXiv.2101.01354, 2021.
        
    \bibitem{GP21} U.K. Gupta, D. Pradhan, Borodin-Kostochka's conjecture on ($P_5, C_4$)-free graphs, J. Appl. Math. Comp., 65 (2021) 877-884.
	
    \bibitem{K80} A.V. Kostochka, Degree, density, and chromatic number, Metody Diskret. Anal. 35 (1980) 45-70 (in Russian).
    
     \bibitem{LLZ24} K. Lan, F. Liu, Y. Zhou, Borodin-Kostochka's Conjecture on $\{P_2\cup P_3,C_4\}$-Free Graphs, Graph. Comb., 40 (2024) 123.
    
	\bibitem{R99} B. Reed, A strengthening of Brooks' theorem, J. Comb. Theory. Ser. B, 76 (1999) 136-149.
	
	\bibitem{W23} H. Wang, The Borodin-Kostochka conjecture for some $(P_2\cup P_3)$-free graphs, Available at SSRN 4412678, 2023.
	
	
	
	
\end{thebibliography}
\end{document}